\documentclass[12pt]{amsart}
\font\bbbld=msbm10 scaled\magstephalf

\newcommand{\bfH}{\hbox{\bbbld H}}

\newcommand{\bfR}{\hbox{\bbbld R}}

\newcommand{\ve}{{\bf e}}
\newcommand{\vn}{{\bf n}}

\newcommand{\ol}{\overline}
\newcommand{\ul}{\underline}

\newcommand{\p}{\partial}
\newcommand{\f}{\frac}

\newtheorem{theorem}{Theorem}[section]
\newtheorem{lemma}[theorem]{Lemma}

\newtheorem{corollary}[theorem]{Corollary}
\theoremstyle{definition}

\theoremstyle{remark}
\newtheorem{remark}[theorem]{Remark}
\numberwithin{equation}{section}

\begin{document}
\setlength{\baselineskip}{1.2\baselineskip}

\title[Entire Spacelike Hypersurfaces]
{Entire Spacelike Hypersurfaces of Prescribed Gauss Curvature in
Minkowski Space}
\author{Bo Guan}
\address{Department of Mathematics, Ohio State University,
         Columbus, OH 43210, USA}
\email{guan@math.osu.edu}
\author{Huai-Yu Jian}
\address{Department of  Mathematical Sciences, Tsinghua University,
Beijing 100084, China}
\email{hjian@math.tsinghua.edu.cn}
\author{Richard M. Schoen}
\address{Department of Mathematics, Stanford University,
          Stanford, CA 94305, USA}
\email{schoen@math.stanford.edu}
\thanks{Research of the first and third authors was supported in part
by NSF grants.
Research of the second author was supported in part by
the National 973-Project and the Trans-Century Training Programme Foundation
for the Talents from the Ministry of Education}


\maketitle

\section{Introduction}
\label{gj2-I}
\setcounter{equation}{0}

In this paper we are concerned with spacelike convex hypersurfaces
of positive constant (K-hypersurfaces) or prescribed Gauss curvature in
Minkowski space
$\bfR^{n, 1}$ ($n \geq 2$).
Any such hypersurface may be written locally as the graph of a convex
function $x_{n+1} = u (x)$, $x \in \bfR^n$ satisfying the spacelike condition
\begin{equation}
\label{gjs-I10}
 |Du| < 1
\end{equation}
and the Monge-Amp\`ere type equation
\begin{equation}
\label{gjs-I20}
\det D^2 u = \psi (x, u) (1 - |Du|^2)^{\frac{n+2}{2}}
\end{equation}
where $\psi$ is a prescribed positive function (the Gauss curvature).
Our main purpose is to study entire solutions on $\bfR^n$
of (\ref{gjs-I10})-(\ref{gjs-I20}).

For $\psi \equiv 1$ a well known entire solution of
(\ref{gjs-I10})-(\ref{gjs-I20}) is the hyperboloid
\begin{equation}
\label{gjs-I30}
 x_{n+1} = \sqrt{1 + |x|^2}, \;\; x \in \bfR^n
\end{equation}
which gives an isometric embedding of the hyperbolic space $\bfH^n$
into $\bfR^{n,1}$.
Hano and Nomizu~\cite{HN83} were probably the first to observe
the non-uniqueness of isometric embeddings of $\bfH^2$ in $\bfR^{2,1}$
by constructing other (geometrically distinct) entire solutions of
(\ref{gjs-I10})-(\ref{gjs-I20}) for $n=2$ (and $\psi \equiv 1$)
using methods of ordinary differential equations.
Using the theory of Monge-Amp\`ere equations, A.-M. Li~\cite{Li95} studied
entire spacelike K-hypersurfaces with uniformly bounded principal curvatures,
while the Dirichlet problem for (\ref{gjs-I10})-(\ref{gjs-I20}) in a
bounded domain $\Omega \subset \bfR^n$ was treated by
Delano\"e~\cite{Delanoe90} when $\Omega$ is strictly convex,
and by Guan~\cite{Guan98a} for general (non-convex) $\Omega$.
In this paper we are interested in entire spacelike K-hypersurfaces,
and more generally hypersurfaces of prescribed Gauss curvature,
without a boundedness assumption on principal curvatures.

Our first goal is to classify all entire spacelike K-hypersurfaces with
symmetries, i.e. those invariant under a subgroup of isometries of
$\bfR^{n,1}$, extending the results of Hano-Nomizu~\cite{HN83} to higher
dimensions. We will focus on hypersurfaces which are rotationally symmetric
with respect to a spacelike axis, as a rotationally symmetric entire spacelike
K-hypersurface with other types of axes either does not exist (when the axis
is lightlike) or is congruent to a rescaling of the standard hyperboloid
(\ref{gjs-I30}) (when the axis is timelike).
These surfaces will be constructed in Section~\ref{gjs-S} where we will
study their properties and asymptotic behavior at infinity. As we will see
in Section~\ref{gjs-M}, understanding these surfaces is crucial to our
study of the Minkowski type problem described below.
One of our main results in Section~\ref{gjs-S} states that these symmetric
K-hypersurfaces are complete with respect to the induced metric from
$\bfR^{n,1}$.

For general entire spacelike K-hypersurfaces it is an important question to
understand their asymptotic behavior at infinity. Li~\cite{Li95} proved that
an entire spacelike
K-hypersurface given by a convex solution $u \in C^{\infty} (\bfR^n)$
of (\ref{gjs-I10})-(\ref{gjs-I20}) has uniformly bounded principal curvatures
if and only if $D u (\bfR^n) = B_1 (0)$, the unit ball in $\bfR^n$.
On the other hand, as we will see in Section~\ref{gjs-S} there do exist
entire K-hypersurfaces with unbounded principal curvatures.
As in the case of hypersurfaces with constant mean curvature
which was treated in \cite{Treibergs82} and \cite{CT90},
the asymptotic behavior of an entire spacelike K-hypersurface can be
characterized by its tangent cone at infinity. (See Section~\ref{gjs-T}.)
Finding entire spacelike K-hypersurfaces with prescribed tangent cones at
infinity is more subtle.
A substantial difficulty is due to the fact that
spacelike K-hypersurfaces do not admit {\em a priori}
interior uniform bounds which keep them from becoming null.
To overcome this difficulty we adopt a variational approach, following
an idea from \cite{GS5}, that allows us to introduce an appropriate class of
weak solutions to (\ref{gjs-I20}), called {\em admissible maximal solutions}
which may only satisfy the weakly spacelike condition
\begin{equation}
\label{gjs-I10'}
 |Du| \leq 1.
\end{equation}
The details will be discussed in Section~\ref{gjs-T} where we consider
the existence and regularity of entire weak solutions to (\ref{gjs-I20})
with prescribed tangent cone at infinity.

Another interesting approach to finding entire spacelike hypersurfaces
with prescribed Gauss curvature and tangent cone at infinity
is to consider the Minkowski type problem of prescribing
Gauss curvature as a function (defined on a
domain $\Omega$ in $\bfH^n$, the unit sphere in $\bfR^{n,1}$) of the unit
normal vector of the prospective hypersurface.
This was indeed the approach employed by Li~\cite{Li95}
who considered the case when the function is defined on the whole space
$\bfH^n$ (or equivalently $B_1 (0) \subset \bfR^n$ via the Legendre
transformation), coupled with a smoothness requirement on the asymptotic
behavior at infinity of the prospective solution $\mbox{graph} (u)$
(in terms of $x \cdot D u (x) - u (x)$).  With the aid of the K-hypersurfaces
constructed in Section~\ref{gjs-S}, we
extend Li's result to allow Lipschitz boundary data for $n=2$, which
geometrically seems to be a more natural assumption.
Another challenging problem is to study more general cases of prescribing the
function on only part of $\bfH^n$.
In this paper we are able to treat the case
$\Omega = \bfH^n_+ := \bfH^n \cap \{x_1 > 0\}$.
This part of the work is included in Section~\ref{gjs-M}.
We hope to come back to the problem in future work.

The corresponding questions for spacelike hypersurfaces of constant mean
curvature have received considerably more intensive investigation.
In their remarkable work on the Bernstein theorem for maximal
hypersurfaces which
extends earlier results due to Calabi~\cite{Calabi68} to higher dimensions,
Cheng-Yau~\cite{CY76} proved that entire spacelike hypersurfaces of constant
mean curvature in $\bfR^{n,1}$ are complete (with respect to the induced
metric) and have uniformly bounded principal curvatures.
Subsequently, Treibergs~\cite{Treibergs82} and Choi-Treibergs~\cite{CT90}
studied the asymptotic behavior at infinity of entire spacelike graphs of
constant mean curvature and treated the existence of such hypersurfaces
with prescribed tangent cone at infinity.
In \cite{BS82} Bartnik-Simon dealt with the Dirichlet problem for the
equation of prescribed mean curvature. Our results seem to indicate that
there are significant differences between entire spacelike hypersurfaces
of constant Gauss curvature and those of constant mean curvature.
It is an interesting open question whether an entire spacelike
K-hypersurface must be complete.

{\bf Acknowledgments.}
Part of this work was done while the second author was visiting the
University of Tennessee and he wishes to thank
the Department of Mathematics for the hospitality.

\bigskip

\section{Entire spacelike K-hypersurfaces with $SO(n-1, 1)$ symmetries}
\label{gjs-S}
\setcounter{equation}{0}

In this section we will classify all entire spacelike K-hypersurfaces
which possess a rotational symmetry with respect to a spacelike axis.
Up to rescaling any such hypersurface is congruent in $\bfR^{n,1}$ to the
graph of a convex solution of (\ref{gjs-I10})-(\ref{gjs-I20}) with
$\psi \equiv 1$ of the form
\begin{equation}
\label{gjs-S10}
 u (x) = \sqrt{f (x_1)^2 + |\bar{x}|^2}, \;\;
 \bar{x} = (x_2, \ldots, x_n), \;\; x = (x_1, \bar{x}) \in \bfR^n
\end{equation}
where $f$ is a positive function defined on $\bfR$.
Geometrically the K-hypersurface
$M := \mbox{graph} (u) \subset \bfR^{n,1}$ is invariant under
the isometries
\begin{equation}
\left( \begin{array}{ccc}
         \cosh \theta &         &  \sinh \theta \\
                      & \varPhi_{n-1} &               \\
         \sinh \theta &         &  \cosh \theta
          \end{array}
  \right), \;\; \theta \in \bfR, \; \varPhi_{n-1} \in SO (n-1).
\end{equation}

We first recall some basic local formulas for the geometric quantities of
spacelike hypersurfaces in the Minkowski space $\bfR^{n,1}$
which is $\bfR^{n+1}$ endowed with the Lorentzian metric
\begin{equation}
\label{gjs-P10}
 ds^2 = \sum_{i=1}^n dx_i^2 - dx_{n+1}^2.
\end{equation}
A spacelike hypersurface $M$ in $\bfR^{n,1}$ is a codimension-one submanifold
whose induced metric is Riemannian.
Locally $M$ can be written as a graph
$x_{n+1} = u (x)$, $x \in \bfR^n$,
satisfying the spacelike condition~(\ref{gjs-I10}).
The induced metric and second fundamental form of $M$ are given by
\begin{equation}
\label{gjs-P20}
 g_{ij} = \delta_{ij} - u_{x_i} u_{x_j}
\end{equation}
and, respectively,
\begin{equation}
\label{gjs-P30}
 h_{ij} = \frac{u_{x_i x_j}}{\sqrt{1 - |Du|^2}},
\end{equation}
while the timelike unit normal vector field to $M$ is
\begin{equation}
\label{gjs-P40}
 \nu =\f {(Du, 1)}{\sqrt{1-|Du|^2}},
\end{equation}
where $Du = (u_{x_1}, \cdots , u_{x_n})$ and $D^2 u = \{u_{x_i x_j}\}$
denote the ordinary gradient and Hessian of $u$, respectively.
We will use $\nabla u$ to denote the gradient of $u$ on $M$.
Note that the norm of $\nabla u$ (with respect to the induce metric
on $M$ from $\bfR^{n,1}$) is
\begin{equation}
\label{gjs-P50}
|\nabla u| \equiv \sqrt{g^{ij} u_{x_i} u_{x_j}} = \frac{|Du|}{\sqrt{1-|Du|^2}}
\end{equation}
where
\begin{equation}
\label{gjs-P60}
 g^{ij} = \delta_{ij} + \frac {u_{x_i} u_{x_j}}{1-|Du|^2}
\end{equation}
is the inverse matrix of $\{g_{ij}\}$.
The Gauss-Kronecker curvature, which is the product of the
principal curvatures (i.e. the eigenvalues of the second fundamental
form with respect to the metric of $M$),
and the mean curvature of $M$ are given by
\begin{equation}
\label{gjs-P70}
 K_M = \frac{\det D^2 u}{(1 - |Du|^2)^{\frac{n+2}{2}}}
\end{equation}
and, respectively
\begin{equation}
\label{gjs-P75}
 H_M = \frac{1}{n} \mbox{div} \Big(\frac{Du}{\sqrt{1-|Du|^2}}\Big).
\end{equation}
Thus equation (\ref{gjs-I20})
locally describes hypersurfaces with prescribed Gauss-Kronecker
curvature $\psi$.

Now assume that $u$ is of the form (\ref{gjs-S10}).
One calculates
\begin{equation}
\label{gjs-S15}
u_{x_1} = \frac{f f'}{u}; \;\;
u_{x_i} =  \frac{x_i}{u}, \;\; 2 \leq i \leq n,
\end{equation}
and
\begin{equation}
\label{gjs-S20}
 1 - |Du|^2 = \frac{f^2 (1 - {f'}^2)}{u^2}.
\end{equation}
Thus $u$ is spacelike if and only if
\begin{equation}
\label{gjs-S25}
|f'| < 1 \;\; \mbox{on $\bfR$}.
\end{equation}
By (\ref{gjs-P50}) and (\ref{gjs-S20}) we have
\begin{equation}
\label{gjs-S27}
 \frac{|\nabla u|}{u} \leq \frac{1}{u \sqrt{1-|Du|^2}}
                             =  \frac{1}{f \sqrt{1 - {f'}^2}}.
\end{equation}
Next,
\begin{equation}
\label{gjs-S30}
\begin{aligned}
  u_{x_1 x_1} & = \frac{f f'' + {f'}^2}{u} - \frac{f^2 {f'}^2}{u^3}
           = \frac{f f'' + {f'}^2 -1}{u} + \frac{g_{11}}{u}, \\
  u_{x_1 x_j} & = - \frac{f f' x_j}{u^3} = \frac{g_{1j}}{u},
               \; 2 \leq j \leq n, \\
  u_{x_i x_j} & = \frac{1}{u} \Big(\delta_{ij} - \frac{x_i x_j}{u^2}\Big)
           = \frac{g_{ij}}{u},  \;\; 2 \leq i, j \leq n
\end{aligned}
\end{equation}
and therefore,
\[ \det D^2 u = \frac{f^3 f''}{u^{n+2}}.  \]
The Gauss curvature of the spacelike hypersurface $M$ in $\bfR^{n,1}$ is
thus given by
\begin{equation}
\label{gjs-S35}
K_M  = \frac{f''}{f^{n-1} (1 - {f'}^2)^{\frac{n+2}{2}}}
\end{equation}
while, by (\ref{gjs-S20}) and (\ref{gjs-S30}), the principal curvatures are
\begin{equation}
\label{gjs-S36}
\kappa_1 = \frac{f''}{(1 - {f'}^2)^{\frac{3}{2}}}, \;\;
\kappa_2 = \ldots = \kappa_n = \frac{1}{f (1 - {f'}^2)^{\frac{1}{2}}}.
\end{equation}
Consequently, if $K_M \equiv 1$
then
\begin{equation}
\label{gjs-S40}
 f'' = f^{n-1} (1 - {f'}^2)^{\frac{n+2}{2}}.
\end{equation}
Integrating (\ref{gjs-S40}) we obtain
\begin{equation}
\label{gjs-S50}
(1-{f'}^2)^{-n/2} - f^n = (1-b^2)^{-n/2} - a^n \equiv c
\end{equation}
where
\begin{equation}
\label{gjs-S60}
 a= f (0), \;\; b = f' (0).
\end{equation}

We summarize some of our observations in the following.

\begin{lemma}
\label{gjs-lemma-S10}
Let  $a > 0$, $|b| < 1$ and $c = (1-b^2)^{-n/2} - a^n$. The following
results hold:

{\bf (a)}
The (unique) solution $f$ to (\ref{gjs-S40}) and (\ref{gjs-S60})
exists on the entire $\bfR$ and satisfies (\ref{gjs-S25}).

{\bf (b)}
If $b \geq 0$ then
\begin{equation}
\label{gjs-S80}
\lim_{t \rightarrow + \infty} f' (t) = 1 \;\; \mbox{and} \;\;
 \lim_{t \rightarrow + \infty} f (t)^2 (1 - f' (t)^2) = 1.
\end{equation}

{\bf (c)}
If $c \leq 1$ then $f > 0$ and $f'' > 0$ on $\bfR$.

{\bf (d)} If $c > 1$ then $f$ changes signs on $\bfR$.

{\bf (e)} Suppose $g$ is another solution of (\ref{gjs-S50}) satisfying
$g (0) > 0$ and $|g' (0)| < 1$. Then either $g \equiv (1 - c)^{1/n}$,
which is possible only when $c < 1$, or there exists $t_0 \in \bfR$
such that $g (t) = f (\alpha t + t_0)$ where $\alpha = 1$ or $-1$.
\end{lemma}

\begin{proof}
 Suppose $f' (t_0) = 1$ for some $t_0 \in \bfR$. We may assume
$t_0 > 0$ and $0 \leq f' < 1$ in $[0, t_0)$. Then
\[ f (t) = f (0) + \int_0^t f' (t) dt < a + t_0,
   \;\; \forall \; 0 \leq t < t_0. \]
However, by (\ref{gjs-S50}),
\[ \lim_{t \rightarrow t_0^-} f (t) = +\infty. \]
This contradiction shows that $|f'| < 1$ wherever the solution exists.
By the theory of ordinary differential equations we see the solution
extends to the entire $\bfR$. This proves (a).

If $b \geq 0$ then from (\ref{gjs-S40}) we see $f'' (t) > 0$ and
$f' (t) > 0$ on $t > 0$. It follows that
\[ \lim_{t \rightarrow + \infty} f (t) = + \infty. \]
By (\ref{gjs-S50}) this implies (\ref{gjs-S80}) and (b) is proved.

From (\ref{gjs-S40}) we see  $f'' > 0$ if $f > 0$ while
$f^n \geq 1 - c$ by (\ref{gjs-S50}).
Now suppose $c = 1$ and $f (t_0) = 0$ for some $t_0 \in \bfR$.
Then $f' (t_0) = 0$ and therefore $f \equiv 0$ by the uniqueness
of solution. This contradicts the fact that $f(0) = a > 0$,
proving (c).

Suppose that $c > 1$ and $f \geq 0$ on $\bfR$. Then $|f'| \geq
(1-c^{-2/n})^{1/2} \equiv \tilde{c} > 0$ on $\bfR$ by
(\ref{gjs-S50}). Without loss of generality, let us assume $f'
\geq \tilde{c}$ on $\bfR$. Then

\[ f (t) = f (0) + \int_0^t f' (t) dt \leq a + \tilde{c} t,
           \;\; \forall \; t \leq 0. \]
Letting $t \rightarrow -\infty$ we reach a contradiction, which implies
(d).

Finally, to prove (e) we observe that if $g$ is not constant then
it also satisfies (\ref{gjs-S40}). From the proof of (b) we see
that $g$ is unbounded above on $\bfR$. There exist therefore $t_1,
t_2 \in \bfR$ such that $f (t_1) = g (t_2)$ and hence $|f' (t_1)|
= |g' (t_2)|$ by (\ref{gjs-S50}). The function $\tilde{f} (t) = f (\alpha (t
- t_2) + t_1)$ where
\[ \alpha = \begin{cases}  1, & \; \mbox{if} \; f' (t_1) = g' (t_2) \\
                          -1, & \; \mbox{if} \; f' (t_1) = - g' (t_2) \neq 0
            \end{cases} \]
then satisfies (\ref{gjs-S40}) and
\[ \tilde{f} (t_2) = g (t_2), \;\; \tilde{f}' (t_2) = g' (t_2). \]
By the uniqueness of solutions we have $\tilde{f} = g$.
The proof is complete.
\end{proof}

By Lemma~\ref{gjs-lemma-S10} when $c > 1$ the corresponding function
$u$ given by (\ref{gjs-S10}) fails to be smooth in $\bfR^n$ while when
$c \leq 1$ the resulting hypersurface is a smooth spacelike
strictly convex entire graph.
Our next lemma enables us to classify these surfaces.

\begin{lemma}
\label{gjs-lemma-S20}
Suppose $a > 0$, $0 \leq b < 1$, $c \equiv (1-b^2)^{-n/2} - a^n \leq 1$
and let $f$ be the solution of (\ref{gjs-S40}) and (\ref{gjs-S60})
on $\bfR$.

{\bf (a)}
If $c = 1$ then $f' > 0$ on $\bfR$ and
\begin{equation}
\label{gjs-S90}
 \lim_{t \rightarrow - \infty} f (t) = 0 \;\; \mbox{and} \;\;
 \lim_{t \rightarrow - \infty} f' (t) = 0.
\end{equation}

{\bf (b)}
If $c < 1$ then there exists $\tau \in \bfR$ such that
$\tilde{f} (t) \equiv f (t + \tau)$ is an even function.
In particular, if $c = 0$
then $\tilde{f} (t) = \sqrt{1 + t^2}$.
\end{lemma}

\begin{proof}
We first consider the case $c = 1$. Suppose $f' (t_0)= 0$ for some
$t_0 \in \bfR$. Then $f (t_0)= 0$ by
(\ref{gjs-S50}) and therefore $f \equiv
0$ by the uniqueness of solution, which is a contradiction. Thus
$f'
> 0$ on the entire $\bfR$. Since $f$ is convex and bounded below
from zero, we have $f' (t) \rightarrow 0$ and hence $f (t)
\rightarrow 0$ by (\ref{gjs-S50}) as $t$ approaches negative
infinity. This proves (a).

Now suppose $c < 1$ and let $h$ be the unique solution of (\ref{gjs-S40})
satisfying $h' (0) = 0$ and $h (0) = (1 - c)^{1/n} > 0$. Then
$h$ is an even function as $h (-t)$ is also a solution of (\ref{gjs-S40})
satisfying the same initial conditions. By Lemma~\ref{gjs-lemma-S10} (e)
we have $h (t) \equiv f (t + \tau)$ for some $\tau \in \bfR$.
\end{proof}

It follows from Lemma~\ref{gjs-lemma-S10} that for each constant $c \leq 1$,
up to a translation and reflection there exists a unique positive
solution $f_c$ of (\ref{gjs-S25}) and (\ref{gjs-S40}) which satisfies
(\ref{gjs-S50}) on $\bfR$.
According to Lemma~\ref{gjs-lemma-S20} we will assume throughout the paper
$f_c$ is even for $c < 1$, and that $f_1$ is chosen so that $f_1 (0) = 1$ and
$f_1' (t) > 0$ for all $t \in \bfR$.
Note that $f_0 (t) = \sqrt{1 + t^2}$.
Let $\mathfrak{H}_c$ denote the graph of
\begin{equation}
\label{gjs-S95}
u_c (x) := \sqrt{f_c (x_1)^2 + |\bar{x}|^2}, \;\; x \in \bfR^n.
\end{equation}
We see that $\mathfrak{H}_c$ is a spacelike entire graph of
constant Gauss curvature one in $\bfR^{n,1}$.
Our main result of this section is the following characterization of
$\mathfrak{H}_c$.

\begin{theorem}
\label{gjs-thm-S10}
{\bf (a)}
For all $c \leq 1$, $\mathfrak{H}_c$
is a complete Riemannian manifold with respect to the induced metric from
$\bfR^{n,1}$.
{\bf (b)}
The principal curvatures of $\mathfrak{H}_c$ are uniformly bounded
for $c < 1$, while $\mathfrak{H}_1$ has unbounded principal curvatures.
{\bf (c)} $D u_c (\bfR^n) = B_1 (0)$ for all $c < 1$ and
$D u_1 (\bfR^n) = B^+_1 (0) := B_1 (0) \cap \{x_1 > 0\}$.
\end{theorem}

\begin{proof}
Note that the principal curvatures are given by (\ref{gjs-S36}).
Part (b) therefore follows from Lemma~\ref{gjs-lemma-S20} and
Lemma~\ref{gjs-lemma-S10} (b), as does part (c) in view of (\ref{gjs-S15}).

To prove part (a) we write $f = f_c$ and $u = u_c$.
Let $\alpha (s) = (x (s), u (s)), s \in [0, L)$ be a geodesic
ray on $\mathfrak{H}_c$ parametrized by arc length such that
$|x (s)| \rightarrow \infty$ as $s \rightarrow L$.
By (\ref{gjs-S27}) we have
\[ \log u (s) - \log u (0) \leq \int_0^s  \frac{|\nabla u|}{u} ds
 \leq \int_0^s \frac{ds}{f \sqrt{1 - {f'}^2}},
                        \;\; \forall \, 0 \leq s < L. \]
If $c < 1$ we see from $f \sqrt{1 - {f'}^2} \geq \sqrt{1 - c}$ that
\[ \log u (s) - \log u (0) \leq \frac{s}{\sqrt{1 - c}},
     \;\; \forall \, s < L. \]
It follows that $L = \infty$ since $u$ is a proper function on $\bfR^n$
in this case.

We now consider case $c = 1$ and assume $f' > 0$.
Suppose there exists some constant $N > 0$ such that
$x_1 (s) \geq -N$ for all $0 \leq s < L$.
We then have  $L = \infty$ as in the previous case ($c < 1$)
since, by Lemma~\ref{gjs-lemma-S20} (a),
$f \sqrt{1 - {f'}^2} \geq c_0 > 0$ for all $0 \leq s < L$ where
$c_0$ is a constant.

Now assume that
\[ \liminf_{s \rightarrow L} x_1 (s) = - \infty. \]
Let $g_{ij}$ be the metric of $\mathfrak{H}_1$. We claim that
\begin{equation}
\label{gjs-S100}
g_{ij} \xi_i \xi_j \geq (1 - (f')^2) \xi_1^2,
        \;\; \forall \; \xi = (\xi_1, \bar{\xi}) \in \bfR^n.
\end{equation}
This follows from the following calculations
\[ g_{11} \xi_1^2 = (1 - (f')^2) \xi_1^2
                    + \frac{(f')^2 |\bar{x}|^2 \xi_1^2}{u^2}  \]
\[ 2 \sum_{i \geq 2} g_{1i} \xi_1 \xi_i
    = - \frac{2 f f' \xi_1}{u^2} \sum_{i \geq 2} x_i \xi_i
    \geq - \frac{(f')^2 |\bar{x}|^2 \xi_1^2}{u^2}
         - \frac{f^2 |\bar{\xi}|^2}{u^2} \]
and
\[  \sum_{i,j \geq 2} g_{ij} \xi_i \xi_j
     = |\bar{\xi}|^2 - \frac{(\bar{x} \cdot \bar{\xi})^2}{u^2}
    \geq \frac{f^2 |\bar{\xi}|^2}{u^2}. \]
Using (\ref{gjs-S100}) we obtain
\[ \begin{aligned}
 s & = \int_0^s
     \Big(g_{ij} \frac{d x_i}{ds} \frac{d x_j}{ds}\Big)^{\frac{1}{2}} ds \\
    & \geq \int_0^s \sqrt{1 - (f')^2} \big|\frac{d x_1}{ds}\big| ds \\
    & \geq - \int_{x_1 (0)}^{x_1 (s)} \sqrt{1 - (f')^2} dx_1 \\
    & \geq - \int_a^{x_1 (s)} \sqrt{1 - (f')^2} dx_1 \\
    & \geq \frac{- x_1 (s) + a}{2}, \;\; \forall \, 0 \leq s < L,
  \end{aligned} \]
where the constant $a \leq x_1 (0)$ is chosen to satisfy
$f' (t) \leq \frac{1}{\sqrt{2}}$ for $t \leq a$.
Letting $s \rightarrow L$ we obtain $L = \infty$.
\end{proof}


\begin{remark}
When $c < 1$ part (a) of Theorem~\ref{gjs-thm-S10} also follows from a
result of Li~\cite{Li95}
as the principal curvatures of $\mathfrak{H}_c$ are bounded.
\end{remark}

\begin{remark}
Up to rescaling any entire spacelike K-hypersurface $M$ in $\bfR^{n,1}$
which is roataitonally symmetric about a spacelike line is
congruent to $\mathfrak{H}_c$ for some $c < 1$ if the principal curvatures
of $M$ are uniformly bounded, and to $\mathfrak{H}_1$ otherwise.
\end{remark}

These K-hypersurfaces 
will be used to construct barrier
function in our study of the Minkowski type problem in
Section~\ref{gjs-M}.
For this purpose we need to know more accurate asymptotic behavior at
infinity of these hypersurfaces. The rest of this section is devoted to
this topic.
Our main tool is the following comparison result for solutions
of (\ref{gjs-S40}).
For a solution  $f$ of (\ref{gjs-S25}), (\ref{gjs-S40}) we denote
$C_f \equiv (1-{f'}^2)^{-n/2} - f^n \leq 1$.

\begin{lemma}
\label{gjs-lemma-S30}
Let $f$ and $g$ be positive solutions of
(\ref{gjs-S25}), (\ref{gjs-S40}) with $C_f < C_g \leq 1$.
Then
{\bf (a)}
$|f' (t)| < |g' (t)|$  wherever $f (t) < g (t)$; and
{\bf (b)}
if $f' (t_0) =  g' (t_0)$ for some $t_0 \in \bfR$
then $f (t) - g (t) \geq f (t_0) - g (t_0) > 0$ for all $t \in \bfR$.
Moreover, $f' (t) > g' (t)$ for all $t > t_0$ and
$f' (t) < g' (t)$ for all $t < t_0$.
\end{lemma}

\begin{proof}
Clearly (a) follows from equation (\ref{gjs-S50}). To prove (b) let
$h = f - g$. Since $C_f < C_g$ we have $h > 0$ by (\ref{gjs-S50}) and,
therefore, $h'' > 0$ by (\ref{gjs-S40}) whenever $h' = 0$.
Consequently, $h$ attains a positive local minimum at any critical point.
This implies that $h$ can have at most one critical point;
(b) is thus proved.
\end{proof}

\begin{corollary}
\label{gjs-cor-S10}
{\bf (a)}
If $c < 0$ or $c = 1$ then
\[ \sqrt{1 + t^2} < f_c (t) < \sqrt{1 + (t + \tau_c)^2}
   \;\;\; \forall \, t > 0  \]
where $\tau_c = \sqrt{(f_c (0))^2 - 1} = \sqrt{(1 - c)^{2/n} - 1}$ for $c < 0$,
and $\tau_1 = \sqrt{2^{2/n} - 1}$.
(Recall that $f_0 (t) = \sqrt{1 + t^2}$.)

{\bf (b)} $0 < c < 1$ then
\[ f_c (t) < \sqrt{1 + t^2} < f_c (t + \tau_c),
   \;\;\; \forall \, t > 0  \]
where $\tau_c > 0$ satisfies $f_c (\tau_c) = f_0 (0) = 1$.
\end{corollary}

\begin{proof}
These are consequences of Lemma~\ref{gjs-lemma-S30}~(b) (applied to
$f_c$ and $f_0$; recall that $f_0 = \sqrt{1 + t^2}$) and the uniqueness
of solutions to the boundary value problems of equation (\ref{gjs-S40}).
\end{proof}

By Lemma~\ref{gjs-lemma-S30} and Corollary~\ref{gjs-cor-S10},
$f_0 (t) - f_c (t)$ is monotone and bounded for $t > 0$. Consequently,
the limit
\[ \lambda_c \equiv \lim_{t \to +\infty} (f_0(t) - f_c (t)) \]
exists for all $c \leq 1$. Note that $\lambda_c < f_0(0) - f_c (0) < 0$ for
$c < 0$, $\lambda_c > f_0 (0) - f_c (0) > 0$ for $c > 0$,
and $\lambda_1 < 0$.

\begin{theorem}
\label{gjs-thm-S20}
For any $c \leq 1$
\begin{equation}
\label{gjs-S300}
  \lim_{t \rightarrow + \infty} (t f_c' (t) - f_c (t)) = \lambda_c,
\end{equation}
while
\begin{equation}
\label{gjs-S350}
  \lim_{t \rightarrow - \infty} (t f_1' (t) - f_1 (t)) = 0.
\end{equation}
\end{theorem}

\begin{proof}
Let $F_c (t) = t f_c' (t) - f_c (t)$.
By the convexity of $f_c$, $F_c' (t) = t f_c'' (t) > 0$ for $t > 0$
and $F_c' (t) = t f_c'' (t) < 0$ for $t < 0$.

Let us first prove
\begin{equation}
  A \equiv \lim_{t \rightarrow - \infty} F_1 (t) = 0.
\end{equation}
The limit exists since $F_1 (t) < 0$ and $F_1' (t) < 0$ for $t < 0$.
Suppose $A < 0$. Since $F_1 (t) < A$ for $t < 0$ and $f_1 (t) \to 0$ as
$t \to -\infty$, there exists $T < 0$ such that
\[  t f_1' (t) <A + f_1 (t) < 0, \ \ \forall \, t \leq T. \]
Thus
\[ \frac{f_1'(t)}{A + f_1 (t)}\leq \frac{1}{t}, \ \ \forall \, t \leq T \]
and
\[ \ln |f_1 (T) + A| - \ln |f_1 (t) + A| \leq \ln |T| - \ln |t|,
   \ \ \forall \, t \leq T. \]
Letting $t \to -\infty$ we obtain a contradiction
\[ \ln |f_1 (T) + A| - \ln |A| = - \infty. \]
This proves (\ref{gjs-S350}).

We next prove (\ref{gjs-S300}) for $c < 0$; the proof for $0 < c \leq 1$ is
similar and will be omitted. In the rest of this proof let $c < 0$ be fixed.
For any fixed $N \geq 0$ there is unique $S_N > 0$ and $T_N > 0$ such that
$f_c' (N) = f_0' (N + S_N)$ and $f_c (N) = f_0 (N + T_N)$.
We have
\begin{equation}
\label{gjs-S330}
 f_0 (t + S_N) + f_c (N) - f_0 (N + S_N)) < f_c (t) < f_0 (t + T_N),
     \;\; \forall \, t > N
\end{equation}
and
\begin{equation}
\label{gjs-S320}
f_0' (t) < f_c' (t) < f_0' (t + T_N),  \;\; \forall \, t > N.
\end{equation}
by Lemma~\ref{gjs-lemma-S30} ((a) for the second inequality in
(\ref{gjs-S320}) and (b) for the first ones in (\ref{gjs-S330})
and (\ref{gjs-S320}). Note that $f_c' (0) = f_0' (0)$.)
Consequently,
\[ \begin{aligned}
   F_c (t) < & t f_0' (t + T_N) - (f_0 (t + S_N) + f_c (N) - f_0 (N + S_N)) \\
           < & F_0 (t + T_N) + f_0 (t + T_N) - T_N f_0' (t + T_N) \\
             & \;\;  -
f_0 (t + S_N) + f_0 (N + S_N) - f_c (N) \\
           < & F_0 (t + T_N)
               + f_0 (t) - f_0 (t + S_N) + f_0 (N + S_N) - f_c (N),
\;\;\;\; \forall \, t > N
    \end{aligned} \]
since $f_0 (t + T_N) - T_N f_0' (t + T_N) < f_0 (t)$ by the convexity of $f_0$.
Thus $\lim F_c (t)$ exists as $t \rightarrow +\infty$ and
\begin{equation}
\label{gjs-S360}
\lim_{t \rightarrow + \infty} F_c (t)
                  \leq f_0 (N + S_N) - f_c (N) - S_N
\end{equation}
as
\[ \lim_{t \rightarrow + \infty} F_0 (t) = 0 \]

and
\[  \lim_{t \rightarrow + \infty} (f_0 (t + S_N) - f_0 (t)) = S_N.  \]
On the other hand,
from (\ref{gjs-S330}) and (\ref{gjs-S320}) we have
\[ F_c (t) > t f_0' (t) - f_0 (t + T_N)
        = F_0 (t) + f_0 (t) - f_0 (t + T_N), \;\;\; \forall \, t > N. \]
It follows that
\begin{equation}
\label{gjs-S370}
\lim_{t \rightarrow + \infty} F_c (t)
     \geq \lim_{t \rightarrow + \infty} (f_0 (t) - f_0 (t + T_N)) = - T_N.
\end{equation}
Note that
\[  \lim_{N \rightarrow + \infty} (f_0 (N + S_N) - f_0 (N) - S_N) = 0  \]
and
\[ 
\lim_{N \rightarrow + \infty} T_N
       = \lim_{N \rightarrow + \infty} (f_0 (N + T_N) - f_0 (N))
       = \lim_{N \rightarrow + \infty} (f_c (N) - f_0 (N)) = - \lambda_c. \]

Letting $N$ approach infinity, from (\ref{gjs-S360}) and (\ref{gjs-S370})
we obtain (\ref{gjs-S300}).
\end{proof}

\begin{corollary}
\label{gjs-cor-S20}
Let $\tilde{f}_1 (t) = f_1 (t + \lambda_1)$. Then
\[ \lim_{|t| \to \infty} (t \tilde{f}_1' (t) - \tilde{f}_1 (t)) = 0. \]
\end{corollary}

\begin{corollary}
\label{gjs-cor-S30}
Let $u_c^*$ be the Legendre transform
of $u$ defined by
\[ u_c^* (y) = \sup \{x \cdot y - u (x): x \in \bfR^n\},
     \;\; y \in  D u_c (\bfR^n).\]
Then
\begin{equation}
\label{gjs-S380}
u_c^* (y) = \begin{cases}
  \lambda_c |y_1|,
     \;\; \mbox{for $y = (y_1, \bar{y}) \in \partial B_1 (0)$, if $c < 1$} \\
  \lambda_c y_1,
     \;\; \mbox{for $y = (y_1, \bar{y}) \in \partial B_1^+ (0)$, if $c = 1$}

    \end{cases}
\end{equation}
where $B_1 (0)$ is the unit ball in $\bfR^n$, and
$B_1^+ (0) = B_1 (0) \cap \{y_1 > 0\}$.
\end{corollary}

\begin{proof}
For any $y \in \Omega_c \equiv D u_c (\bfR^n)$, by (\ref{gjs-S15})
\[ u_c^* (y) = x \cdot D u_c (x) - u_c (x)
     = \frac{f_c (x_1) (x_1 f_c' (x_1) - f_c (x_1))}{u_c (x)}, \]
where $x = (x_1, \bar{x}) \in \bfR^n$ is uniquely given by
$D u_c (x) = y$. Letting $y$ approach an arbitrarily fixed point on
$\partial \Omega_c$ we obtain (\ref{gjs-S380}) from Theorem~\ref{gjs-thm-S20}
and (\ref{gjs-S15}).
\end{proof}

This proves to be useful in Section~\ref{gjs-M}
where we will also need the following lemma

\begin{lemma}
\label{gjs-lemma-S50}
$\lambda_c \to - \infty$ as $c \to - \infty$ and
$\lambda_c \to + \infty$ as $c \to 1^-$.
\end{lemma}

\begin{proof}
The first case is obvious since
$\lambda_c < f_0 (0) - f_c (0) = 1 - (1-c)^{1/n}$ for $c < 0$.
Next, for any fixed $N > 0$ there exists $c_N \in (0, 1)$ such that
\[ f_c (0) = (1-c)^{1/n} < f_1 (-2 N), \;\;\; \forall \, c_N < c < 1. \]
By Lemma~\ref{gjs-lemma-S30} (a),
\[ f_c (t) < f_1 (t-2 N), \;\;\; \forall \, t > 0, \, c_N < c < 1. \]
In particular,
\[ f_c (N) < f_1 (- N), \;\;\; \forall \, c_N < c < 1. \]
It follows that
\[ \lambda_c > f_0 (N) - f_c (N) > f_0 (N) - f_1 (- N),
                                 \;\;\; \forall \, c_N < c < 1 \]
since $f_0 (t) - f_c (t)$ is increasing for $t > 0$ when $c > 0$.
Letting $c \to 1^-$ and then $N \to +\infty$, we prove the second case.
\end{proof}

\bigskip

\section{The tangent cone at infinity}
\label{gjs-T}
\setcounter{equation}{0}

In this section we first characterize the tangent cones for
entire spacelike convex hypersurfaces in Minkowski space with bounded
Gauss curvature.
We then will consider the problem of finding such K-hypersurfaces of with
a prescribed tangent cone at infinity.
Let $u$ be an entire convex solution of
(\ref{gjs-I10})-(\ref{gjs-I20}) with
$0 < \psi_1 \leq \psi \leq \psi_2$ on $\bfR^n$ where $\psi_1$, $\psi_2$
are constant. Consider
\[ u_r (x) := \frac{u(r x)}{r}, \;\; x \in \bfR^n, \; r > 0, \]
\begin{equation}
\label{gjs-T10}
V_u (x):=\lim_{r \to 0} u_r (x), \;\; x \in \bfR^n.
\end{equation}
Following \cite{CT90} and \cite{Treibergs82} we call $V_u$ the {\em blowdown}
of $u$ at infinity.
Note that, by (\ref{gjs-I10}) and the convexity of $u$,
$V_u$ is well-defined and convex on $\bfR^n$,
\begin{equation}
\label{gjs-T20}
V_u (\lambda x) = \lambda V_u (x), \;\; \forall \; x \in \bfR^n, \; \lambda > 0
\end{equation}
and
\begin{equation}
\label{gjs-T30}
|V_u(x)-V_u(y)|\leq |x-y|, \;\; \forall \; x, y \in \bfR^n.
\end{equation}
Moreover, $V_u$ satisfies the {\em null condition}, that is

\begin{lemma}
\label{gjs-lemma-K10}
For any $x \in \bfR^n$ there exists
$y \in \bfR^n$, $y \neq x$, such that
\begin{equation}
\label{gjs-T40}
|V_u (x) - V_u (y)| = |x-y|.
\end{equation}
\end{lemma}

\begin{proof}
Suppose this is not true. Then there exists $x_0\in \bfR^n$
and $\delta >0$ such that
\[ V_u (x) \leq V_u(x_0) + 1 - 2 \delta,
  \;\; \forall \, x \in \partial B_1 (x_0) \]
where $B_1 (x_0)$ is the unit ball in $\bfR^n$ centered at $x_0$.
By the convexity of $u$ we have
\[ \frac{d}{dr} (u_r (x) - u_r (0)) \leq 0, \;\; \forall x \in \bfR^n. \]
Thus the limit in (\ref{gjs-T10}) is uniform on compact sets by Dini's
Theorem.
Consequently, we can find $r_0>0$ such that
\begin{equation}
u_r (x) \leq V_u (x_0) + 1 - \delta, \;\; \forall \, x \in \partial B_1 (x_0)
\end{equation}
for all $r > r_0$.
It therefore follows from the maximum principle that
\[ u_r (x) \leq W (x; r)
 := V_u (x_0) + ((\psi_1^{1/n} r)^{-2} + |x - x_0|^2)^{\f {1}{2}} - \delta,
   \;\; \forall \, x \in B_1 (x_0) \]
as both $u_r$ and $W (\cdot; r)$ are spacelike
in $B_1 (x_0)$ and
\[ \det D^2 u_r (x)   =  r^n \det D^2 u (r x)
                   \geq  r^n \psi_1 (1 - |D u|^2)^{\frac{n+2}{2}},
                      \; \;\; x \in B_1 (x_0) \]
while
\[ \det D^2 W (x; r) = r^n \psi_1 (1 - |D W (x;r)|^2)^{\frac{n+2}{2}},
                     \;\; x \in B_1 (x_0). \]
Letting $r \to \infty$ we obtain
\[ V_u (x_0) \leq V_u(x_0) - \delta, \]
which is a contradiction.
\end{proof}

Recall that the set of subdifferentials of a convex function $v$
at a point $x_0 \in \bfR^n $ is defined as
\[ T_v (x_0) := \{\alpha \in \bfR^n:
     v(x) \geq v (x_0) + \alpha \cdot (x - x_0), \; \forall x \in \bfR^n\}. \]
Obviously, $T_v (x_0)$ is a closed convex set and equals
$Dv (x_0)$ if $v$ is differentiable at $x_0$.
We call $\overline {T_{V_u}(\bfR^n)}$ the {\em tangent cone at infinity} of
graph $u$. Using Lemma~\ref{gjs-lemma-K10}
one can show as in \cite{CT90} that
\begin{equation}
\label{gjs-T50}
\ol{T_{V_u}(\bfR^n)} = T_{V_u}(0) = \ol{Du(\bfR^n)} \subseteq \ol{B_1 (0)}
\end{equation}
and
\begin{equation}
\label{gjs-T60}
V_u (y) = |y|, \;\; \forall \, y \in \ol{Du(\bfR^n)}.
\end{equation}
This last identity can be seen as follows.
By definition
\[ V_u (y) \geq V_u (0) + y \cdot y = |y|^2,
\;\; \forall \, y \in T_{V_u}(0) \]
since $V_u (0) = 0$. In particular, from (\ref{gjs-T30}) we have
\[ V_u (y) = 1, \;\; \forall \, y \in  T_{V_u}(0) \cap \partial B_1 (0) \]
By (\ref{gjs-T20}), we therefore obtain (\ref{gjs-T60}).
The following lemma can also be shown as in \cite{CT90}.

\begin{lemma}
\label{gjs-lemma-T20}
$T_{V_u}(0)$ is the convex hull of $T_{V_u}(0) \cap \partial B_1 (0)$.
In particular, $T_{V_u}(0)$ has no interior strictly extremal points.
Moreover,
\[ V_u (x) = \sup \{\alpha \cdot x :
      \alpha \in  T_{V_u}(0) \cap \partial B_1 (0)\},
    \;\; x \in \bfR^n. \]
\end{lemma}

It is a natural question to find entire K-hypersurfaces with a given
tangent cone. In order to treat this problem we introduce a class of
weak solutions to (\ref{gjs-I20}) and discuss their basic properties.

For a domain $\Omega \subseteq \bfR^n$ and a nonnegative function $\psi$
defined on $\Omega \times \bfR$, let $\mathcal{A}[\psi, \Omega]$ denote the
collection of weakly spacelike, locally convex subsolutions (in the viscosity
sense) of (\ref{gjs-I20}) in $C^0(\ol{\Omega})$.
We call $u \in \mathcal{A}[\psi, \Omega]$ an {\em admissible maximal solution}
of (\ref{gjs-I20}) in $\Omega$ if
\begin{equation}
\label{gjs-T100}
 \int_{\Omega'} \sqrt{1 - |Du|^2} dx
       \geq \int_{\Omega'} \sqrt{1 - |Dv|^2} dx
\end{equation}
for any bounded subdomain $\Omega'$ of $\Omega$
and $v \in \mathcal{A}[\psi, \Omega']$ with $u = v$ on $\partial \Omega$.
Note that (\ref{gjs-T100}) means geometrically that the volume of
graph$_{\Omega'} (u)$ is greater than or equal to that of
graph$_{\Omega'} (v)$. Thus the graph of an admissible maximal solution
is a volume maximizer in $\mathcal{A}[\psi, \Omega]$.

\begin{lemma}
\label{gjs-lemma-M10}
Let $u \in \mathcal{A}[\psi, \Omega]$ be an admissible maximal solution of
(\ref{gjs-I20}).
If $u$ is spacelike 
in a subdomain
$\Omega' \subseteq \Omega$, then it is a viscosity solution in $\Omega'$.
In particular, if $u \in C^2 (\Omega')$ then it is a classical solution, and
is  locally strictly convex if $\psi > 0$.
\end{lemma}

\begin{proof}
We first assume that $\Omega'$ is smooth and bounded,
$\psi \in C^{\infty} (\ol{\Omega'} \times \bfR)$, $\psi > 0$, and
$u \in C^2 (\ol{\Omega'})$.
Using $u$ as a subsolution, we can apply a theorem in \cite{Guan98a} to obtain
a spacelike locally strict convex solution $v \in  C^{\infty} (\ol{\Omega'})$
of (\ref{gjs-I20}) satisfying $v \geq u$ in $\ol{\Omega'}$ and
$v = u$ on $\partial \ol{\Omega'}$. By Lemma~\ref{gjs-lemma-M20} (below)
we have
\[  \int_{\Omega'} \sqrt{1 - |Du|^2} dx
       \leq \int_{\Omega'} \sqrt{1 - |Dv|^2} dx.  \]
Replacing $u$ by $v$ on $\Omega'$, we obtain a function
$\tilde{u} \in \mathcal{A}[\psi, \Omega]$.
By the definition of admissible maximal solutions we see that the equality
holds and therefore $v = u$ in $\Omega'$.
By an approximation argument we prove the lemma in the general case.
\end{proof}

\begin{lemma}
\label{gjs-lemma-M20}
Let $u_1, u_2 \in C^{0,1} (\Omega) \cap C^0 (\ol{\Omega})$ be spacelike and
satisfy $u_1 \geq u_2$ in $\ol{\Omega}$ and $u_1 = u_2$ on $\partial \Omega$.
Suppose $u_1$ is convex, or more generally,
the spacelike graph of $u_1$ in $\bfR^{n,1}$ has nonnegative
generalized mean curvature almost everywhere, that is
\[ \mbox{div} \Big(\frac{D u_1}{\sqrt{1 - |D u_1|^2}}\Big) \geq 0 \;\; a.e. \]
Then
\[ \int_{\Omega} \sqrt{1 - |D u_1|^2} dx
       \geq \int_{\Omega} \sqrt{1 - |D u_2|^2} dx.  \]
The equality holds if and only if $u_1 = u_2$ in $\Omega$.
\end{lemma}

\begin{proof}
Let $S_i$ denote the graph of $u_i$ in $\bfR^{n+1}$ over $\Omega$ and
\[ \nu_i = \frac{(- D u_i (x), 1)}{\sqrt{1 + |D u_i (x)|^2}} \]
the (Euclidean) upward unit normal vector field to $S_i$, $i =1, 2$.
Consider the vector filed
\[ N (x, z) = \frac{(D u_1 (x), 1)}{\sqrt{1 - |D u_1 (x)|^2}},
\;\; (x, z) \in R \]
where
\[ R := \{(x,z) \in \bfR^{n+1}:
            u_2 (x) < z < u_1 (x), \; x \in \Omega\} \]
is the region in $\bfR^{n+1}$ bounded by $S_1$ and $S_2$.
We have
\[ \mbox{div} N (x, z)
   = \mbox{div} \Big(\frac{D u_1}{\sqrt{1 - |D u_1|^2}}\Big) \geq 0
  \;\;\; a.e. /; \mbox{in $R$}. \]
Consequently by the divergence theorem
\[ \begin{aligned}
0 \leq \int_R \mbox{div} N  dv
   & = \int_{S_1} N \cdot \nu_1 d \sigma
       - \int_{S_2} N \cdot \nu_2 d \sigma \\
   & =  \int_{\Omega} \sqrt{1 - |D u_1|^2} dx
       - \int_{\Omega} \frac{1 - D u_1 \cdot D u_2}{\sqrt{1 - |D u_1|^2}} dx \\
   & \leq \int_{\Omega} \sqrt{1 - |D u_1|^2} dx
         - \int_{\Omega} \sqrt{1 - |D u_2|^2} dx.
   \end{aligned} \]
The last inequality follows from
\[ (1 - D u_1 \cdot D u_2)^2 \geq (1 - |D u_1|^2) (1 - |D u_2|^2). \]
Obviously, all the equalities hold if and only if $u_1 = u_2$ in $\Omega$.
\end{proof}

We now state our existence result of this section.

\begin{theorem}
\label{gjs-thm-T10}
Let $E$ be a subset of $\partial B_1 (0)$ which is not contained in
any hyperplane in $\bfR^n$.
Then there exists a convex admissible
maximal solution $u \in C^{0,1} (\bfR^n)$ to (\ref{gjs-I20}) with
$\psi \equiv 1$ satisfying
\begin{equation}
\label{gjs-T70}
 \overline {Du (\bfR^n)} = \Gamma (E),
\end{equation}
where $\Gamma (E)$ denotes the convex hull of $E$, and
\begin{equation}
\label{gjs-T80}
V_u (x) = V_E := \sup_{\alpha \in E} \alpha \cdot x, \;\; x \in \bfR^n.
\end{equation}
\end{theorem}

\begin{proof}
By a theorem of Choi-Treibergs~\cite{CT90} there exists a
spacelike entire graph $x_{n+1} = v (x)$, $v \in C^{\infty} (\bfR^n)$,
of mean curvature one whose tangent cone is $\Gamma (E)$. Moreover,
$v$ is strictly convex and satisfies $v \geq V_v = V_E$  on $\bfR^n$.

For each integer $k \geq 1$, by a theorem of Delano\"e~\cite{Delanoe90}
there exists a unique spacelike strictly convex solution
$u_k \in C^{\infty} (\ol{B_k (0)})$ to the Dirichlet problem
\[ \begin{aligned}
 \det D^2 u & = (1 - |Du|^2)^{\frac{n+2}{2}} \;\; \mbox{in $\ol{B_k (0)}$} \\
          u & = v \;\; \mbox{on $\partial B_k (0)$}.
   \end{aligned} \]
Since $|Du_k| \leq 1$ and $|D V_E| = 1$ where $D V_E$ exists, by the maximum
principle we have $V_E \leq u_k \leq v$ on $\ol{B_k (0)}$ for all $k$.
Moreover, there exists a subsequence $u_{k_j}$
and a weakly spacelike convex function $u \in C^{0,1} (\bfR^n)$ such that
$u_{k_j}$ converges to $u$ in $C^{0,1} (\ol{\Omega})$ for any bounded domain
$\Omega$ in $\bfR^n$. It follows from Lemma~\ref{gjs-lemma-M20} and the
comparison principle that $u$ is an admissible maximal solution to
(\ref{gjs-I20}). Note that $V_E \leq u \leq v$.
From $V_v = V_E$ we obtain (\ref{gjs-T80}) and therefore (\ref{gjs-T70})
by (\ref{gjs-T50}).
\end{proof}

\section{The Minkowski type problem}
\label{gjs-M}
\setcounter{equation}{0}

In this section we consider the Minkowski type problem which provides a
natural approach to the problem of finding entire spacelike hypersurfaces
of prescribed Gauss curvature.
Let $M = \mbox{graph} (u)$ be a smooth spacelike strictly convex
hypersurface. Then the Gauss map
\[ \nu: M \to \bfH^n \subset \bfR^{n,1}, \;\;
    \nu (x, u (x)) = \frac{(Du, 1)}{(1-|Du|^2)^{1/2}} \]
is a diffeomorphism from $M$ onto its image in $\bfH^n$.
On the other hand, $\bfH^n$ can be identified with the unit ball $B_1 (0)$
in $\bfR^n$ by the diffeomorphism
\[ \pi: \bfH^n \to B_1 (0), \;\;
  \pi (\xi, \xi_{n+1})
   = \frac{\xi}{\xi_{n+1}}. \]
For convenience we will also call $\vn := \pi \circ \nu$ the Gauss map.
It is immediately seen that
\[ \vn (x, u (x)) = D u (x), \;\; \forall \, x \in \bfR^n. \]

Thus geometric quantities of $M$ can be viewed as defined via the Gauss
map on its image $\Omega := \vn (M) \subseteq B_1 (0)$.
Naturally one can consider the Minkowski type problem: given a domain
$\Omega \subseteq B_1 (0)$ and a function $\eta > 0$ on $\Omega$,
find an entire spacelike strictly convex hypersurface $M = \mbox{graph} (u)$
whose Gauss map image is $\Omega$ and Gauss curvature at $\vn^{-1} (y)$
is given by $\eta (y)$ for $y \in \Omega$ where
$\vn^{-1}: \Omega \to M$ is the inverse Gauss map.

As $\Omega$ has nonempty boundary (in $\bfR^n$), one needs to impose
certain boundary conditions in order to describe the
asymptotic behavior of the hypersurface at infinity. To formulate such a
boundary value problem, we consider the support function of the graph
of $u$ given by the Lorentz inner product $\langle X,\nu\rangle=
(x\cdot Du-u)/\sqrt{1-|du|^2}$. The expression $x \cdot Du (x) - u (x)$, $x \in \bfR^n$
leads us to consider the Legendre transform of $u$
\[ u^* (y) = \sup_{x \in \bfR^n} (x \cdot y - u (x)),
   \;\; y \in \Omega. \]
where $\Omega = Du (\bfR^n) \subseteq B_1 (0)$.
It is well known
that $u^*$ is strictly convex and that for $y \in \Omega$
\[ u^* (y) = x \cdot y - u (x), \;\; Du^* (y) = x \]
and

\[ D^2 u^* (y) = (D^2 u (x))^{-1} \]
where $x \in \bfR^n$ is uniquely determined by $Du (x) = y$.
By (\ref{gjs-I20}) we see that $u^*$ should satisfy the Monge-Amper\`e equation
\begin{equation}
\label{gjs-T110}
\det D^2 v (y) = \frac{1}{\eta (y) (1-|y|^2)^{\frac{n+2}{2}}},
   \;\; \forall \, y \in \Omega
\end{equation}
where $\eta(y)=\psi(x)$.

Conversely, given a convex domain $\Omega \subseteq B_1 (0)$ and
$\eta \in C^{\infty} (\Omega)$, $\eta > 0$, if there exists a strictly
convex solution $v \in C^{\infty} (\Omega)$ of (\ref{gjs-T110}) such that
\begin{equation}
\label{gjs-T120}
Dv (\Omega) = \bfR^n,
\end{equation}
then its Legendre transform $u = v^*$ is a smooth
spacelike strictly convex solution of (\ref{gjs-I20}) defined on $\bfR^n$
with $\psi (x) = \eta (y)$, where $y$ is given by $D v (y) = x$,
for all $x \in \bfR^n$.
According to Li~\cite{Li95}, the resulting hypersurface $M = \mbox{graph} (u)$
has uniformly bounded principal curvatures if and only if $\Omega = B_1 (0)$.

Li~\cite{Li95} treated the Dirichlet problem in $\Omega = B_1 (0)$ for
(\ref{gjs-T110})-(\ref{gjs-T120}) with smooth boundary data. From the
geometric point of view, it would be natural to consider Lipschitz
boundary data, as well as general subdomains of $B_1 (0)$. Analytically,
this is a challenging problem as one has to construct more sophisticated
barrier functions to prove that (\ref{gjs-T120}) is satisfied. (In \cite{Li95}
the barriers are constructed from the function $\sqrt{1 - |y|^2}$ which is
the Legendre transform of the hyperboloid (\ref{gjs-I30}).)
Our main results of this section extend the theorem of
Li~\cite{Li95} to allow Lipschitz boundary data in dimension $n=2$
(Theorem~\ref{gjs-thm-T30}),
and to the case $\Omega = B^+_1 (0)$ (Theorem~\ref{gjs-thm-T20}) for all $n$.
This is achieved with the aid of the rotationally symmetric K-hypersurfaces
$\mathfrak{H}_c$ constructed in Section~\ref{gjs-S}.
We first consider the case $\Omega = B^+_1 (0)$:
write $\partial \Omega = \partial_+ \Omega \cup \partial_0 \Omega$
where $\partial_+ \Omega = \partial \Omega \cap \{y_1 > 0\}$
and $\partial_0 \Omega = \partial \Omega \cap \{y_1 = 0\}$.

\begin{theorem}
\label{gjs-thm-T20}
Let $\Omega = B_1^+ (0)$ and
$\varphi \in C^0 (\partial \Omega) \cap C^{\infty} (\ol{\partial_+ \Omega})$,
$\eta \in C^{\infty} (\Omega) \cap C^0 (\ol{\Omega})$, $\eta > 0$.
Suppose in addition that
\begin{equation}
\label{gjs-T130}
\mbox{$\varphi$ is affine on $\partial_0 \Omega$}.
\end{equation}
Then there exists a unique strictly convex solution
$v \in C^{\infty} (\Omega) \cap C^0 (\ol{\Omega})$ of (\ref{gjs-T110})
which satisfies (\ref{gjs-T120}) and
the Dirichlet condition
\begin{equation}
\label{gjs-T140}
v = \varphi \;\; \mbox{on $\partial \Omega$}.
\end{equation}
\end{theorem}

\begin{proof}
For convenience we write $\psi = 1/\eta$ and will still use $\varphi$ to
denote its harmonic extension to $\ol{\Omega}$.
Note that $\varphi \in C^{\infty} (\Omega \cup \partial_+ \Omega)$.
Let $\Omega_1 \subset \cdots \subset \Omega_k \subset \cdots \subset \Omega$
be a sequence of smooth strictly convex domains such that
\begin{equation}
 \bigcup_{i=1}^{\infty} \Omega_k = \Omega.
\end{equation}
Let $\varepsilon_k \to 1$ be a strictly increasing sequence.
By \cite{CNS1} there exists a unique strictly convex solution
$v_k \in C^{\infty} (\ol{\Omega_k})$ to the Dirichlet problem
\begin{equation}
\label{gjs-T150}
\left\{ \begin{aligned}
\det D^2 v_k & = \psi (1- \varepsilon_k |y|^2)^{-\f{n+2}{2}}
                 \;\; \mbox{in $\ol{\Omega_k}$} \\
        v_k & = \varphi \;\; \mbox{on $\partial \Omega_k$}.
\end{aligned} \right.
\end{equation}
By the maximum principle
\begin{equation}
\label{gjs-T160}
\varphi \geq v_k > v_{k+1} \geq \ul{v} \;\;
          \mbox{in $\ol{\Omega_k}$}, \; \forall \, k \geq 1.
\end{equation}
where
\[ \ul{v} (y) = \ul{\varphi} - \bar{\psi}^{\frac{1}{n}} \sqrt{1 - |y|^2},
                \;\; y \in \ol{B_1}, \]
\[ \ul{\varphi} = \min_{\partial \Omega} \varphi, \;\;
   \bar{\psi} = \max_{\ol{\Omega}} \psi, \]
since $\ul{v}$ is a subsolution of (\ref{gjs-T150}) for each $k \geq 1$, i.e.
\begin{equation}
\det D^2 \ul{v} = \bar{\psi} (1- |y|^2)^{-\f{n+2}{2}}
           \geq \psi (1- \varepsilon_k |y|^2)^{-\f{n+2}{2}}
                 \;\; \mbox{in $\ol{\Omega_k}$}
\end{equation}
and $\ul{v} \leq \varphi$ on $\partial \Omega_k$.
From (\ref{gjs-T160}) we obtain by the convexity of $v_k$ a uniform bound
on any compact subset of $\Omega$ for $|D v_k|$ independent of $k$.
It follows that $v_k$ converges uniformly on any compact set in $\Omega$
to the convex function $v \in C^0 (\Omega)$ given by
\[ v (y) = \lim_{k \to \infty} v_k (y), \;\; y \in \Omega. \]

Next, for an arbitrarily fixed point $\hat{y} \in \partial \Omega$
by subtracting an affine function we may assume
$\varphi (\hat{y}) = 0$ and $D \varphi (\hat{y}) = 0$.
Since
$\varphi \in C^0 (\partial \Omega) \cap C^{\infty} (\ol{\partial_+ \Omega})$
and $\varphi$ is affine on $\partial_0 \Omega$
we can choose $A > 0$ sufficiently large depending on
$|D \varphi|_{\ol{\partial_+ \Omega}}$ such that
\begin{equation}
\label{gjs-T170}
 - A l (y) \leq \varphi (y) \leq A l (y)
          \;\; \forall \, y \in \partial \Omega
\end{equation}
where $l (y) = 1 - \hat{y} \cdot y$ if $\hat{y} \in \partial_+ \Omega$,
$l (y) = y_1$ if $\hat{y} \in \partial_0 \Omega$.
By the maximum principle we have as in (\ref{gjs-T160}) that
\begin{equation}
\label{gjs-T180}
\varphi (y) \geq v_k (y)
                \geq  \bar{\psi}^{\frac{1}{n}} u_1^* (y) - A l (y),
           \;\; \forall \, y \in \Omega_k, \;\; \forall \, k \geq 1.
\end{equation}
Here, with a slight abuse of notation, $u_1^*$ is the Legendre transform of
the function $\tilde{u}_1 (x) := (\tilde{f}_1 (x_1)^2 + |\bar{x}|^2)^{1/2}$
where $\tilde{f_1} (t) = f (t + \lambda_1)$ as in Corollary~\ref{gjs-cor-S20},
noting that $u_1^* \in C^0 (\ol{B_1^+}) \cap C^{\infty} (B_1^+)$
satisfies
\[ \det D^2 u_1^* = (1 - |y|^2)^{-\frac{n+2}{2}} \;\; \mbox{in $B_1^+$} \]
and $u_1^* = 0$ on $\partial \Omega$ by Corollary~\ref{gjs-cor-S20}.
Letting $k \to \infty$ we obtain from (\ref{gjs-T180}) that
\begin{equation}
\label{gjs-T190}
\lim_{y \to \hat{y}} v (y) = \varphi (\hat{y}),
   \;\; \forall \, \hat{y} \in \partial \Omega
\end{equation}
since $\varphi (\hat{y}) =
\bar{\psi}^{\frac{1}{n}} u_1^* (\hat{y}) - A l (\hat{y}) = 0$.

This proves $v \in C^0 (\ol{\Omega})$ with $v = \varphi$ on $\partial \Omega$.
We next want to prove $v \in C^{\infty} (\Omega)$.
Note that $v$ is a convex viscosity solution of (\ref{gjs-T110}) in $\Omega$.
Let $y_0$ be any interior point in $\Omega$ and $P$ a supporting plane of
$\Sigma_v : = \mbox{graph} (v)$ at $(y_0, v (y_0))$.
We claim that $P \cap \Sigma_v$ contains a single point $(y_0, v (y_0))$.
For otherwise, by a theorem of Caffarelli~\cite{Caffarelli90a},
$P \cap \Sigma_v$ would contain a segment from $(y_0, v (y_0))$
to a boundary point $(\hat{y}, v (\hat{y}))$ for some
$\hat{y} \in \partial \Omega$, which would imply
\begin{equation}
\label{gjs-T200}
 \lim_{t \to 0^+} \frac{v (\hat{y} + t \ve) - v (\hat{y})}{t}
   = \frac{v (y_0) - v (\hat{y})}{|y_0 - \hat{y}|} > - \infty
\end{equation}
where $\ve$ is the unit vector pointing from $y_0$ to $\hat{y}$.
However, by the maximum principle and the second inequality in
(\ref{gjs-T170}) which we may still assume to hold,
\begin{equation}
\label{gjs-T205}
v (y) \leq A l (y) + \ul{\psi}^{\frac{1}{n}} u_1^* (y),
       \;\; \forall \, y \in \ol{\Omega},
\end{equation}
where
\[ \ul{\psi} = \min_{\ol{\Omega}} \psi > 0. \]
It follows that
\begin{equation}
\label{gjs-T210}
 \lim_{t \to 0^+} \frac{v (\hat{y} + t \ve) - v (\hat{y})}{t}
   \leq A \ve \cdot D l + \ul{\psi}^{\frac{1}{n}}
     \lim_{t \to 0^+} \frac{u_1^* (\hat{y} + t \ve) - u_1^* (\hat{y})}{t}
 = - \infty
\end{equation}
since $|D u_1^*| = \infty$ on $\partial \Omega$.
This contradicts (\ref{gjs-T200}), proving our claim.
By Caffarelli's theorems~\cite{Caffarelli90a}, \cite{Caffarelli90b} and
the Evans-Krylov regularity theory $v$ is a smooth strictly convex solution
of (\ref{gjs-T110}) in $\Omega$.
Moreover, from (\ref{gjs-T210}) which holds for any interior point
$y_0 \in \Omega$ and $\hat{y} \in \partial \Omega$, we see $v$ satisfies
(\ref{gjs-T120}).
\end{proof}

\begin{remark}
The resulting entire spacelike hypersurface $M = \mbox{graph} (v^*)$ must
have unbounded principal curvatures.
\end{remark}

\begin{remark}
Assumption (\ref{gjs-T130}) is also necessary when $n = 2$. In
general ($n \geq 2$) it is necessary to assume $\varphi$ to be convex but
not strictly convex at each interior point of $\partial_0 \Omega$.
This is because if $\varphi$ is smooth and strictly convex at a point
$\hat{y} \in B_1 (0) \cap \{y_1 = 0\}$ then the solution is at least
of class $C^{0,1}$ up to boundary near $\hat{y}$ by the boundary regularity
of Monge-Amp\`ere equations.
In particular, (\ref{gjs-T120}) can not hold at $\hat{y}$.
\end{remark}

\begin{remark}
Concerning problem (\ref{gjs-T110})-(\ref{gjs-T120}) in a general subdomain
$\Omega$ of $B_1 (0)$, Lemma~\ref{gjs-lemma-T20} gives a necessary condition
on $\Omega$ for its solvability.
In particular, when $n = 2$ it implies $\Omega$ has to be either $B_1 (0)$
or $B_1 (0) \cap \{a \cdot y > c\}$ for some $a \in \bfR^n$, $|a| = 1$ and
$- 1 < c < 1$. In all dimensions ($n \geq 2$) this latter case can be
reduced to $\Omega = B_1^+ (0)$.
\end{remark}

As we mentioned above, our second main theorem of this section concerns
the Minkowski type problem with Lipschitz Dirichlet boundary data.

\begin{theorem}
\label{gjs-thm-T30}
Let $n=2$, $\Omega = B_1 (0) \subset \bfR^2$,
$\eta \in C^{\infty} (\Omega) \cap C^0 (\ol{\Omega})$, $\eta > 0$,
and $\varphi \in C^{0,1} (\partial \Omega)$.
Then there exists a unique strictly convex solution
$v \in C^{\infty} (\Omega) \cap C^0 (\ol{\Omega})$ of (\ref{gjs-T110})
which satisfies (\ref{gjs-T120}) and (\ref{gjs-T140}).
Consequently, there exists a smooth complete entire spacelike strictly convex

hypersurface $M$ with Gauss curvature
\[ K_M (\vn^{-1} (y)) = \eta (y), \;\; \forall \, y \in B_1 (0) \]
where $\vn^{-1}: B_1 (0) \to M$ is its inverse Gauss map.
\end{theorem}

\begin{proof}
We modify the proof of Theorem~\ref{gjs-thm-T20}. First by approximation
(solving (\ref{gjs-T150}) for $\Omega_k = B_1 (0)$ for all $k \geq 1$)
we obtain a convex viscosity solution $v \in C^0 (\Omega)$ of (\ref{gjs-T110}).
To proceed let $\hat{y} \in \partial \Omega$. We may assume $\hat{y} = (0, 1)$
and $\varphi (\hat{y}) = 0$.
Since $\varphi \in C^{0,1} (\partial \Omega)$, by Corollary~\ref{gjs-cor-S30}
and Lemma~\ref{gjs-lemma-S50} there exists
$c_1 < 0$, $0 < c_2 < 1$ and $A > 0$ (independent of $\hat{y}$) such that
\begin{equation}
\label{gjs-T250}
\bar{\psi}^{\frac{1}{n}} u_{c_1}^* - A (1-y_2)
    \leq \varphi \leq \ul{\psi}^{\frac{1}{n}} u_{c_2}^* + A (1-y_2)
\;\; \mbox{on $\partial \Omega$}.
\end{equation}
Applying the maximum principle to the approximation we obtain
\begin{equation}
\label{gjs-T260}
\bar{\psi}^{\frac{1}{n}} u_{c_1}^* - A (1-y_2)
    \leq v \leq \ul{\psi}^{\frac{1}{n}} u_{c_2}^* + A (1-y_2)
\;\; \mbox{in $\ol{\Omega}$}.
\end{equation}
This proves $v \in C^0 (\ol{\Omega})$ and $v = \varphi$ on $\partial \Omega$.

Finally, using the second inequality in (\ref{gjs-T260}) (in place of
(\ref{gjs-T205})) we can prove $v \in C^{\infty} (\Omega)$ and satisfies
(\ref{gjs-T120}) as in the proof of Theorem~\ref{gjs-thm-T20}.
\end{proof}

It would be interesting to extend Theorem~\ref{gjs-thm-T30} to higher
dimensions.

\end{document}